# Kings and serfs in oriented graphs

[1]S. Pirzada and [2]N. A. Shah

[1,2]Department of Mathematics, University of Kashmir, Srinagar, India
[1]Email: sdpirzada@yahoo.co.in

**Abstract.** In this paper, we extend the concept of kings and serfs in tournaments to that of weak kings and weak serfs in oriented graphs. We obtain various results on the existence of weak kings(weak serfs) in oriented graphs, and show the existence of n-oriented graphs containing exactly k weak kings(weak serfs), $1 \le k \le n$. Also, we give the existence of n-oriented graphs containing exactly k weak kings and exactly s weak serfs such that b weak kings from k are also weak serfs.

A tournament is a digraph in which every pair of distinct vertices is joined by exactly one arc. If there is a path of length 1 or 2 from a vertex u to a vertex v, then v is said to be reachable from u within two steps. A vertex u in a tournament is called a king if every other vertex v in T is reachable within two steps from u. A vertex u is called a serf if u is reachable within two steps from every other vertex v in T. Maurer [9] introduced the dual terms king and serf in a delightful exposition of a tournament model for dominance in flocks of chickens. Landau [7] proved that every vertex of maximum score in a tournament is a king, and in [2] it can be seen that a tournament without transmitters contains at least three kings. The following result is due to Maurer [9] and an inductive proof of it is given by Reid [15].

**Theorem 1.** For all integers $n \geq k \geq 1$ there exists an n-tournament with exactly k kings with the following exceptions: k=2 with n arbitrary, and n=k=4 (in which case no such n-tournament exists).

Moon [10] proved that no tournament has exactly two kings. Maurer [9] asked to determine all 4-tuples (n,k,s,b) for which there exists an n-tournament with exactly k kings and s serfs such that b of the kings are also serfs. Such a tournament is called a (n,k,s,b)-tournament. The following characterization of such 4-tuples is given by Reid [14].

**Theorem 2.** Suppose that $n \geq k \geq s \geq b \geq 0$ and $n > 0$. There exists a (n,k,s,b)-tournament if and only if the following conditions hold.
(1) $n \geq k + s - b$, (2) $s \neq 2$ and $k \neq 2$,
(3) either $n = k = s = b \neq 4$ or $n > k$ and $s > b$,
(4) (n,k,s,b) is none of (n,4,3,2),(5,4,1,0),or (7,6,3,2).

Maurer [9] also asked to determine those n-tournaments T which are contained in a tournament whose kings are exactly the vertices of T. The following characterization is due to Maurer [9] and the proof can be found in [15].

**Theorem 3.** A nontrivial n-tournament T is contained in a tournament whose kings are the vertices of T if and only if T contains no transmitter.

Further Reid [15] determines the smallest order m so that there exists an m-tournament W which contains the given n-tournament T as a subtournament so that every vertex of W is a king. Reid [15] also obtains bounds in a similar problem in which the kings of W are exactly the vertices of T, and Lee and Chang [8] improved such bounds.

An r-king ($r \geq 2$) of a digraph D is a vertex of D from which any other vertex of D can be reached in less or equal to r steps. Gutin [1] proved that each k-partite ($k \geq 2$) tournament with at most one transmitter contains a 4-king and that there exist infinitely many multipartite tournaments without 3-king. The same was rediscovered independently by Petrovic and Thomassen [11]. Koh and Tan [3,4,5,6] obtained several results on kings in multipartite tournaments. Petrovic [12] presents a variety of results concerning r-kings (r=2,3,4) of bipartite tournaments.

An oriented graph is a digraph with no symmetric pairs of directed arcs and without loops. Define $s_{v_i}$ (or $s(v_i)) = n - 1 + d^+_{v_i} - d^-_{v_i}$, the score of a vertex $v_i$ in an oriented graph D, where $d^+_{v_i}$ and $d^-_{v_i}$ are the outdegree and indegree, respectively, of $v_i$ and n is the number of vertices in D. The score sequence of an oriented graph is formed by listing the vertex scores in non-decreasing order.

The following terms and notations can be found in [13]. For any two vertices u and v in an oriented graph D, we have one of the following possibilities.

(i) An arc directed from u to v, denoted by u(1–0)v.
(ii) An arc directed from v to u, denoted by u(0–1)v.
(iii) There is no arc from u to v and there is no arc from v to u, and is denoted by u(0–0)v.

If $d^*(v)$ is the number of those vertices u in D which have v(0–0)u, then $d^+(v) + d^-(v) + d^*(v) = n - 1$. Therefore, $s(v) = 2d^+(v) + d^*(v)$. This implies that each vertex u with v(1–0)u contributes two to the score of v, and each vertex u with v(0–0)u contributes one to the score of v. Since the number of arcs and non-arcs in an oriented graph of order n is $\binom{n}{2}$, and each v(0–0)u contributes two (one each at u and v) to scores, therefore the sum total of all the scores is $2\binom{n}{2}$. An oriented graph can be interpreted as the result of a round-robin competition in which ties (draws) are allowed, that is, the participants play each other once, with an arc from u to v if and only if u defeats v. A player receives two points for each win , and one point for each tie . With this scoring system , player v obtains a total of s(v) points .

A triple in an oriented graph is an induced oriented subgraph with three vertices. For any three vertices u, v and w, the triples of the form u(1–0) v(1–0)w(1–0)u, or u(1–0) v(1–0)w(0–0)u, or u(0–0)v(1–0)w(1–0)u, or u(1–0)v(0–0)w(1–0)u are said to be intransitive, while as the triples of the form u(1–0)v(1–0)w(0–1)u, or u(0–1)v(1–0)w(1–0)u, or u(1–0)v(0–1)w(1–0)u, or u (1–0)v(0–1)w(0–0)u , or u(0–1)v(0–0)w(1–0)u, or u(0–0)v(1–0)w(0–1)u, or u(1–0)v(0–0)w(0–1)u, or u(0–0)v(0–1)w(1–0)u, or u(0–1)v(1–0)w(0–0)u, or u(1–0)v(0–0)w(0–0)u, or u(0–1)v(0–0)w(0–0)u, or u(0–0)v(1–0)w(0–0)u, or u(0–0)v(0–1)w(0–0)u, or u(0–0)v(0–0)w(1–0)u, or u(0–0)v(0–0)w(0–1)u are said to be transitive. An oriented graph is said to be transitive if all its triples are transitive. The converse $D^{/}$ of an oriented graph D is obtained by reversing each arc of D.

Let u and v be vertices in an oriented graph D such that either u(1-0)v, or u(0-0)v, or u(1-0)w(1-0)v, or u(1-0)w(0-0)v, or u(0-0)w(1-0)v for some vertex w in D. Then v is said to be weakly reachable within two steps from u. If either u(1-0)v, or u(1-0)w(1-0)v for some w in D, then v is reachable within two steps from u.

A vertex u in an oriented graph D is called a weak king if every other vertex v in D is weakly reachable within two steps from u. A vertex u is called a king if every other vertex v in D is reachable within two steps from u.

A vertex u in an oriented graph D is called a weak serf if u is weakly reachable within two steps from every other vertex in D, and a vertex u in D is called a serf if u is reachable within two steps from every other vertex v in D.

We observe in the null oriented graph, every vertex is a weak king as well as a weak serf. Therefore, the oriented graphs considered in this paper are not null unless otherwise stated.

Also, we note that there exist n-oriented graphs with exactly k kings for all integers $n \geq k \geq 1$, only with the exception $n = k = 4$. Theorem 1 guarantees the existence of complete n-oriented graphs(tournaments) with exactly k kings for all integers $n \geq k \geq 1$, with the exception: k=2 and n arbitrary, and $n = k = 4$. An oriented graph D with exactly two kings is constructed as follows. Let V={$v_1,v_2,…,v_n$} be vertex set of D, and arcs defined as $v_1(1-0)v_i$, for i = 2,4,…,n ; $v_1(0-1)v_3$; $v_2(1–0)v_3$ ; $v_2(1-0)v_i$, for $4 \leq i \leq n$ ; and for all other $i \neq j$, $v_i(0-0)v_j$. The vertices $v_1$ and $v_3$ are kings in D.

There are no complete oriented graphs(tournaments) with 4 vertices and exactly 4 kings, and therefore there do not exist incomplete oriented graphs with 4 vertices and exactly 4 kings. This is because the incomplete oriented graph can be obtained from the complete one by removing its one or more arcs which further destroys the reachability.

Here we study weak kings and weak serfs, and in this direction we have the following result.

**Theorem 4.** If u is a vertex with maximum score in an oriented graph D, then every other vertex v in D is weakly reachable within two steps from u.

**Proof.** Let u be a vertex of maximum score in D, and let X, Y, and Z be respectively the set of vertices x, y, and z such that u(1-0)x, u(0-0)y, and u(0-1)z. Let $|X| = n_1, |Y| = n_2$, and $|Z| = n_3$. Clearly, s(u)=$2n_1+n_2$. If $n_3=0$, the result is trivial. So, assume $n_3 \neq 0$. Claim that each z Є Z is weakly reachable within two steps from u. If not, let $z_0$ be a vertex in Z not weakly reachable within two steps from u. Then for each x Є X and each y Є Y, $z_0(1-0)x$, and $z_0(1-0)y$ or $z_0(0-0)y$. In case $z_0(1-0)x$ and $z_0(1-0)y$ for each x Є X and each y Є Y, then s($z_0$)$\geq 2+2n_1+2n_2$ = s(u)+$n_2$+2 > s(u), which contradicts the choice of u. If $z_0(1-0)x$ and $z_0(0-0)y$ for each x Є X and each y Є Y, then s($z_0$)$\geq 2+2n_1+n_2$ = s(u)+2 > s(u), again contradicting the choice of u. This establishes the claim, and hence the proof is complete. ■

The next result is a direct consequence of Theorem 4.

**Theorem 5.** A vertex of maximum score in an oriented graph is a weak king.

In tournaments, the vertex of maximum score is always a king. This need not be true in oriented graphs. To see this, consider the oriented graph D shown in Figure 1. The scores of vertices $v_1$, $v_2$, $v_3$ and $v_4$ are respectively 2,3,3 and 4. Clearly, $v_4$ is a vertex of maximum score but is not a king as $v_1$ is not reachable within two steps from $v_4$. However, $v_4$ is a weak king.

Now, consider the oriented graph $D^*$ with vertices $u_1$, $u_2$, $u_3$, $u_4$ and $u_5$, and arcs defined by $u_1(1-0)u_2$, $u_2(1-0)u_i$ for i = 3,4,5 and $u_i(0-0)u_j$ for all other i ≠ j. Clearly, $s(u_1) = 5$, $s(u_2) = 6$, $s(u_3) = 3$, $s(u_4) = 3$, and $s(u_5) = 3$. Evidently $u_1$ is a king in $D^*$ while as the vertex $u_2$ of maximum score is not a king.

A vertex v in an oriented graph is a transmitter if indegree of v is zero, that is v is a transmitter if $d^-(v)=0$.

**Theorem 6.** An oriented graph D without transmitters has atleast three weak kings.
**Proof.** Let be a vertex of maximum score in the oriented graph D. Clearly, by Theorem 5, u is a weak king. As D has no transmitters, there is atleast one vertex v such that v(1-0)u. Let V be the set of these vertices v, and let $v_1$ be the vertex of maximum score in V. Let X, Y and Z respectively be the set of vertices x, y and z , other than u, with $v_1(1-0)x$, $v_1(0-0)y$ and $v_1(0-1)z$. Assume $|X|=n_1, |Y|=n_2$, and $|Z|=n_3$ so that $s(v_1)=2n_1+n_2+2$. Claim that all vertices of Z are weakly reachable within two steps from $v_1$. If not, let $z_0$ be a vertex which is not weakly reachable within two steps from $v_1$. Then $z_0(1-0)u$, and (a) $z_0(1-0)x$ and (b) $z_0(1-0)y$ or $z_0(0-0)y$ for each x Є X and each y Є Y.

If for each x in X and each y in Y, $z_0(1\text{-}0)x$ and $z_0(1\text{-}0)y$, then $s(z_0) \geq 2n_1+2n_2+4 = s(v_1)+n_2+2 > s(v_1)$. This contradicts the choice of $v_1$. If for each x in X and each y in Y, $z_0(1\text{-}0)x$ and $z_0(0\text{-}0)y$, then $s(z_0) \geq 2n_1+n_2+4 > s(v_1)$, again contradicting the choice of $v_1$. This establishes the claim, and thus $v_1$ is also a weak king.

Now, let W be set of vertices w with $w(1\text{-}0)v_1$ and let $w_1$ be the vertex of maximum score in W. Then by the same argument as above, every other vertex in D can weakly reach within two steps from $w_1$, and so $w_1$ is a weak king. Since D is asymmetric, and in D we have $w_1(1\text{-}0)v_1$ and $v_1(1\text{-}0)u$, necessarily u, $v_1$ and $w_1$ are distinct vertices. Hence, D contains at least three weak kings. ∎

The following result shows that indeed n-oriented graph with exactly k weak kings, for all integers $n \geq k \geq 1$, do exist.

**Theorem 7.** There exist n-oriented graphs with exactly k weak kings for all integers $n \geq k \geq 1$.

**Proof.** Define the arcs of n-oriented graph D with vertex set $U = \{x, y, u_1, u_2, \ldots, u_{n-2}\}$ as

$x(0\text{-}0)y$,

$u_i(1\text{-}0)x$ and $u_i(0\text{-}1)y$ for all $1 \leq i \leq n-2$,

$u_i(0\text{-}0)u_j$ for all $i \neq j$ and $1 \leq i \leq n-2$, $1 \leq j \leq n-2$.

Clearly, x is a weak king as $x(0\text{-}0)y$ and $x(0\text{-}0)y(1\text{-}0)u_i$ for all i, $1 \leq i \leq n-2$. Also, y is a weak king as $y(0\text{-}0)x$ and $y(1\text{-}0)u_i$ for all i, $1 \leq i \leq n-2$. Finally, every $u_i$ is a weak king, since $u_i(0\text{-}0)u_j$ for all $i \neq j$ and $u_i(1\text{-}0)x$ and $u_i(1\text{-}0)x(0\text{-}0)y$. Thus, D contains exactly n weak kings.

If in D, $u_{n-2}(1\text{-}0)u_{n-3}$, then D contains exactly n-1 weak kings, since $u_{n-2}$ is not weakly reachable within two steps from $u_{n-3}$, and so $u_{n-3}$ is not a weak king.

If in D, $u_{n-2}(1\text{-}0)u_{n-3}$ and $u_{n-2}(1\text{-}0)u_{n-4}$, then D contains exactly n-2 weak kings, with $u_{n-3}$ and $u_{n-4}$ being not weak kings, since $u_{n-2}$ is not weakly reachable within two steps neither from $u_{n-3}$ nor from $u_{n-4}$, Continuing in this way, if in D, $u_{n-2}(1-0)u_i$ for all $k-2 \leq i \leq n-3$, then D contains exactly k weak kings.

If D is constructed so that $x(0-0)y$, $x(1-0)u_i$, $y(1-0)u_i$ and $u_i(0-0)u_j$ for all $1 \leq i \leq n-2$, $1 \leq j \leq n-2$ and $i \neq j$, then D contains exactly two weak kings x and y.

Clearly D has exactly one weak king if it is constructed such that $x(0-0)y$, $u_1(1-0)x$, $u_1(1-0)y$ and $u_1(1-0)u_j$ for all $2 \leq i \leq n-2$ with $x(1-0)u_i$ and $y(1-0)u_i$ for all $2 \leq i \leq n-2$. ∎

Figure 2(a) shows 6 vertex oriented graph with exactly 6 weak kings, 2(b) shows 6 vertex oriented graph with exactly 5 weak kings namely x, y, $v_1$, $v_2$ and $v_4$, 2(c) shows 6 vertex oriented graph with exactly 4 weak kings namely x, y, $v_1$ and $v_4$, 2(d) shows 6 vertex oriented graph with exactly 3 weak kings namely x, y and $v_4$, 2(e) shows 6 vertex oriented graph with exactly 2 weak kings namely x and y, and 2(f) shows 4 vertex oriented graph with exactly 4 weak kings and exactly 4 weak serfs.

The directional dual of a weak king is a weak serf, and thus a vertex u is a weak king of an oriented graph D if and only if u is a weak serf of $\overline{D}$, the converse of D. So by duality and Theorem 7, there exists an n-oriented graph with exactly s weak serfs for all integers $n \geq s \geq 1$. If $n = k \geq 1$ in Theorem 7, then every vertex in any such n-oriented graphs is both a weak king and a weak serf. Also, if $n > k \geq 1$ in Theorem 7, the n-oriented graphs described in the proof of Theorem 7 contains vertices which are both weak kings and weak serfs, and also contains vertices which are weak kings but not weak serfs and vice versa. These ideas give rise to the following problem. For what 4-tuples (n, k, s, b) does there exist an n-oriented graph with exactly k weak kings, s weak serfs and that exactly b of the weak kings are also serfs ? Such oriented graphs are called (n, k, s, b)-oriented graphs. Without loss of generality, we assume that $k \geq s$.

We note that if u is a weak king (weak serf) in an oriented graph D, then u is weakly reachable from (can weakly reach) every other weak king (weak serf) within two steps.

We have the following observation.

**Lemma 1.** Let u be a vertex in an oriented graph D.
  (a) If u is not a weak king, then D contains a vertex v such that v(1-0)u, v is not a weak serf, and arc vu is in no intransitive triples of D.
  (b) If u is not a weak serf, then D contains a vertex w such that u(1-0)w, w is not a weak king, and arc uw is in no intransitive triples of D.

**Proof.** (a) Assume u is not a weak king. Then, there exists a vertex v such that v(1-0)u and y is not weakly reachable within two steps from u. So, v is not a weak serf and u(0-0)x(0-0)v, or u(0-1)x(0-1)v, or u(0-1)x(0-0)v, or u(0-0)x(0-1)v for every vertex x in D. Combining with v(1-0)u, for every x in D, the triples are of the form v(1-0)u(0-0)x(0-0)v, or v(1-0)u(0-1)x(0-1)v, or v(1-0)u(0-1)x(0-0)v, or v(1-0)u(0-0)x(0-1)v. Thus, every triple in D which contains arc vu is transitive. Part (b) can be proved similarly.    ∎

**Theorem 8.** If $n > k \geq s \geq 0$, then there exists no (n, k, s, s)-oriented graph.

**Proof.** Assume $n > k \geq s > 0$ and let D be an (n, k, s, s)-oriented graph. Since every weak serf of D is also a weak king, the suboriented graph of non-weak kings, denoted by N, contains no weak serf. As $n > k$, N contains at least one vertex, and let u be a vertex of N with minimum score in D. Since u is not a serf, by Lemma 1(b), there exists a vertex v such that u(1-0)v and arc uv is in no intransitive triple of D. Therefore, v is not a weak king and so v is in N. Since uv is in no intransitive triple, so u(1-0)x(0-1)v(0-1)u, or u(1-0)x(0-0)v(0-1)u, or u(0-0)x(1-0)v(0-1)u, or u(0-0)x(0-0)v(0-1)u for every vertex x in D. Thus in all cases, s(u) > s(v) which contradicts the choice of u, Hence, there exists no (n, k, s, s)-oriented graphs.    ∎

**Theorem 9.** There exist (n, k, s, b)-oriented graphs, $n \geq k \geq s > b \geq 0$ and $n > 0$, $n \geq k+s-b$.

**Proof.** Let $D_1$ be the oriented graph with vertex set $\{x_1, y_1, u_1, u_2, \ldots, u_{k-b-2}\}$ and $x_1(0\text{-}0)y_1$, $u_i(1\text{-}0)x_1$, $u_i(0\text{-}1)y_1$ for all $1 \leq i \leq k-b-2$, and $u_i(0\text{-}0)u_j$ for all $i \neq j$.

Take the oriented graph $D_2$ with vertex set $\{x_2, y_2, v_1, v_2, \ldots, v_{b-2}\}$ and arcs defined as in $D_1$. Let $D_3$ be the oriented graph with vertex set $\{z_1, z_2, \ldots, z_{s-b}\}$ and $z_i(0\text{-}0)z_j$ for all $i$, $j$. Let $D$ be the oriented graph $D_1 \cup D_2 \cup D_3$ (See Figure 3) with

$z_i(1\text{-}0)y_2$ for $1 \leq i \leq s-b$

$z_i(0\text{-}0)x_2$ for $1 \leq i \leq s-b$

$z_i(0\text{-}0)v_j$ for $1 \leq i \leq s-b$, $1 \leq j \leq b-2$

$x_1(1\text{-}0)z_i$, $y_1(1\text{-}0)z_i$ for $1 \leq i \leq s-b$

$u_r(0\text{-}0)z_i$ for $1 \leq r \leq k-b-2$, $1 \leq i \leq s-b$

$x_1(1\text{-}0)y_2$, $y_1(1\text{-}0)y_2$

$v_r(1\text{-}0)y_2$ for $1 \leq r \leq k-b-2$

$x_1(0\text{-}0)x_2$, $y_1(0\text{-}0)x_2$

$v_r(0\text{-}0)v_j$, for $1 \leq r \leq k-b-2$, $1 \leq j \leq b-2$.

Clearly $D$ contains exactly $k$ weak kings and the weak king set is $\{x_1, y_1\} \cup \{u_1, u_2, \ldots, u_{k-b-2}\} \cup \{x_2, y_2\} \cup \{v_1, v_2, \ldots, v_{b-2}\}$. $D$ contains exactly $s$ weak serfs with the weak serf set as $\{x_2, y_2\} \cup \{v_1, v_2, \ldots, v_{b-2}\} \cup \{z_1, z_2, \ldots, z_{s-b}\}$. Also, from these $k$ weak kings, exactly $b$ are weak serfs. The weak king-serf set is $\{x_2, y_2\} \cup \{v_1, v_2, \ldots, v_{b-2}\}$. ∎

**Theorem 10.** An n-oriented graph $D$ is contained in an oriented graph whose weak kings are the vertices of $D$ if $D$ contains no transmitter.

**Proof.** Without loss of generality, assume that $D$ to be connected (strongly or weakly). Let $D$ have no transmitters. Obviously $n \geq 3$. Let $U = \{u_1, u_2, \ldots, u_n\}$ be the vertex set of $D$. Let $D_1$ be an isomorphic, disjoint copy of $D$ with vertex set $V = \{v_1, v_2, \ldots, v_n\}$ where $v_i$ corresponds to $u_i$, $1 \leq i \leq n$. Consider the oriented graph $D_2 = D \cup D_1$ and $u_i(1\text{-}0)v_j$ for $i \neq j$ and $u_i(0\text{-}1)u_j$ for $i = j$, $1 \leq i \leq n$ and $1 \leq j \leq n$. Clearly, each vertex $u_i$ is a weak king in $D_2$, but no $v_i$ is a weak king of $D_2$, as $u_i$ is not a transmitter of $D$, $1 \leq i \leq n$. ∎

The converse of Theorem 10 is not true in general. To see this, consider the 6-vertex oriented graph D shown in Figure 4. Here x, y, $u_1$ and $u_4$ are weak kings while as $u_2$ and $u_3$ are not. But the suboriented graph induced by the vertices x, y, $u_1$ and $u_4$ do contain a transmitter namely y .

We conclude with the following problem.

**Problem 1.** Obtain the necessary and sufficient conditions for an n-oriented graph to contain exactly k weak kings .

**Acknowledgement.** The authors would like to thank the referee for his helpful suggestions.

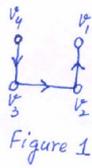

Figure 1

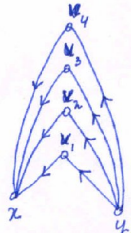

Figure 2(a)

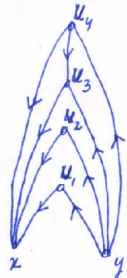

Figure 2(b)

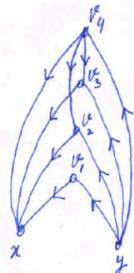

Figure 2(c)

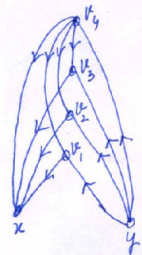

Figure 2(d)

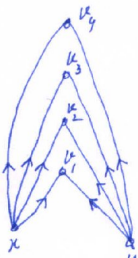

Figure 2(e)

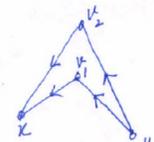

Figure 2(f)

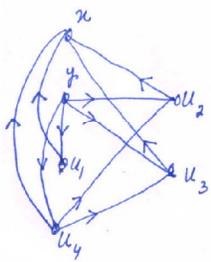

Figure 4

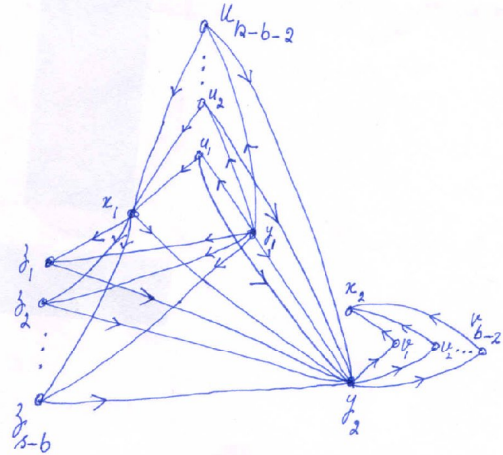

Figure 3

# References

[1]  G. M. Gutin, The radii of n-partite tournaments, Math. Notes 40 (1986) 743-744.

[2]  F. Harary, R. Z. Norman and D. Cartwright, Structural Models : An Introduction to the Theory of Directed graphs, John Wiley and Sons, Inc. New York (1965)294.

[3]  K. M. Koh and B.P.Tan, Kings in multipartite tournaments , Discrete Mathematics, 147(1995)171-183.

[4]  K. M. Koh and B.P.Tan, Number of 4-kings in bipartite tournaments with no 3-kings, Discrete Mathematics, 154(1996) 281-287.

[5]  K. M. Koh and B.P.Tan, The number of kings in multipartite tournaments, Discrete Mathematics, 167(1997) 411-418.

[6]  K. M. Koh and B.P.Tan, The set and number of kings in a multipartite tournament, Bull. Inst. Combin. Appl. 13(1995) 15-22.

[7]  H.G. Landau, On dominance relations and the structure of animal societies: I. Effect of inherent characteristics, Bull. Math. Biophysc. 13(1951) 1-19, II. Some effects of possible social factors, Bull. Math. Biophysc. 13(1951) 245-262, III. The condition for a score structure, Bull. Math. Biophysc. 15(1953) 143-148.

[8]  H.Y.Lee and G.J.Chang, Medians of graphs and kings of tournaments, Taiwanese J.Mathematics, Vol. 1, No. 1 (1997)103-110.

[9]  S.B.Maurer, The king chicken Theorems, Math. Mag. 53(1980) 67-80.

[10]  J.W.Moon, Solution to problem 463, Math. Mag. 35(1962) 189.

[11]  V. Petrovic and C. Thomassen, Kings in k-partite tournaments, Discrete Mathematics, 98 (1991) 237-238.

[12] V. Petrovic, Kings in bipartite tournaments, Discrete Mathematics, 173 (1997)187-196.

[13]  S.Pirzada and U. Samee, Score sequences in oriented graphs, Submitted.

[14]  K.B.Reid, Tournaments with prescribed number of kings and serfs, Congressus Numerantium, Vol. 29 (1980) 809-826.

[15]  K.B.Reid, Every vertex a king, Discrete Mathematics, 38 (1982) 93-98.